# CHARACTERIZATION OF BAYES PROCEDURES FOR MULTIPLE ENDPOINT PROBLEMS AND INADMISSIBILITY OF THE STEP-UP PROCEDURE[1]

By Arthur Cohen and Harold B. Sackrowitz

*Rutgers University*

The problem of multiple endpoint testing for $k$ endpoints is treated as a $2^k$ finite action problem. The loss function chosen is a vector loss function consisting of two components. The two components lead to a vector risk. One component of the vector risk is the false rejection rate (FRR), that is, the expected number of false rejections. The other component is the false acceptance rate (FAR), that is, the expected number of acceptances for which the corresponding null hypothesis is false. This loss function is more stringent than the positive linear combination loss function of Lehmann [*Ann. Math. Statist.* **28** (1957) 1–25] and Cohen and Sackrowitz [*Ann. Statist.* (2005) **33** 126–144] in the sense that the class of admissible rules is larger for this vector risk formulation than for the linear combination risk function. In other words, fewer procedures are inadmissible for the vector risk formulation. The statistical model assumed is that the vector of variables $\mathbf{Z}$ is multivariate normal with mean vector $\boldsymbol{\mu}$ and known intraclass covariance matrix $\Sigma$. The endpoint hypotheses are $H_i:\mu_i=0$ vs $K_i:\mu_i>0$, $i=1,\dots,k$. A characterization of all symmetric Bayes procedures and their limits is obtained. The characterization leads to a complete class theorem. The complete class theorem is used to provide a useful necessary condition for admissibility of a procedure. The main result is that the step-up multiple endpoint procedure is shown to be inadmissible.

**1. Introduction.** Let $\mathbf{Z}$ be a $k\times 1$ random vector which is multivariate normal with mean vector $\boldsymbol{\mu}=(\mu_1,\dots,\mu_k)'$ and known covariance matrix $\Sigma$. Assume $\Sigma$ is intraclass, that is, all variances are equal to $\sigma^2$ and all correlations are equal to $\rho$. Consider the $k$ hypothesis testing problems $H_i:\mu_i=0$

Received February 2003; revised March 2004.
[1]Supported by NSA Grant MDA 904-02-1-0039.
*AMS 2000 subject classifications.* 62C10, 62C15.
*Key words and phrases.* Step-up procedure, intraclass correlation, Bayes procedures, inadmissibility, finite action problem, Schur convexity, complete class.

vs $K_i : \mu_i > 0$, $i = 1, \ldots, k$. This represents the problem of multiple endpoint testing. We view this problem as a $2^k$ finite action problem where we can decide whether to reject or accept each $H_i$ individually. The loss function chosen is a vector loss consisting of two components. The first component is the number of false rejections and the second component is the number of false acceptances. The corresponding vector risk has one component related to the average power of the procedure while the other component is related to the average size of the procedure. This will be made precise in Section 2.

The vector loss function is more stringent than the linear combination loss function used in Lehmann (1957) and Cohen and Sackrowitz (CS) (2005) in the sense that the class of admissible procedures for the vector loss function contains the class of admissible procedures for the linear combination loss function. In other words, any procedure shown to be inadmissible for the vector loss is inadmissible for the linear combination loss.

In this paper we offer a characterization of the class of symmetric (permutation invariant) Bayes procedures. For the normal model, intraclass is the most general case of permutation invariance. The characterization leads to a useful complete class theorem. The complete class theorem yields a useful necessary condition for admissibility of a procedure. Our most important result is that the popular step-up multiple endpoint testing procedure is inadmissible.

Step-up procedures are studied in many places including Hochberg and Tamhane (1987), Hochberg (1988) and Shaffer (1995). Step-up procedures have been studied in connection with procedures that control the false discovery rate (FDR). See Benjamini and Hochberg (1995), Benjamini and Yekutieli (2001) and Sarkar (2002). Six of the eighteen multiple endpoints procedures studied by Dudoit, Shaffer and Boldrick (2003) are step-up procedures.

The inadmissibility result for step-up is of practical importance. Furthermore, the result is somewhat akin to the Stein-type inadmissibility phenomenon in the following sense: The step-up procedure leads to admissible tests for each component individually when $\rho > 0$. See CS (2005). Yet in the finite action problem if the loss function is the vector loss function of this paper or the sum of losses for the component problems, the step-up procedure is inadmissible for $k \geq 2$.

We note that the characterization of Bayes procedures for finite action problems has only been achieved in the past for the case where the action space is a subset of the real line. See, for example, Karlin and Rubin (1956), Brown, Cohen and Strawderman (1976) and Van Houwelingen and Verbeek (1985). Finite action formulations are realistic, practical and important.

Preliminaries and notation will be given in Section 2. The characterization of symmetric Bayes procedures will be given in Section 3. Section 4 contains a description of a complete class, a necessary condition for admissibility, and the result that the step-up procedure is inadmissible. In Section 5, for the case $k=2$ a procedure that is better than step-up is constructed.

**2. Preliminaries.** This $2^k$ finite action problem has actions $\mathbf{a}=(a_1,a_2,\dots,a_k)'$ where $a_i$ equals 0 or 1 for $i=1,\dots,k$. An action where $a_i=1$ means that $H_i$ is rejected, where if $a_i=0$, $H_i$ is accepted. Thus, for example, $\mathbf{a}=(1,\dots,1)'$ means all $H_i$ are rejected. It will be convenient to define

$$\Gamma=\{\mathbf{u}:\mathbf{u}=(u_1,\dots,u_k)',u_i=0 \text{ or } 1, \text{ all } i\}.$$

Note that $\Gamma$ can be used to represent the totality of all actions. However, $\Gamma$ will serve other purposes as well.

Decision rules $\delta(\cdot|\mathbf{z})$ are probability mass functions on $\Gamma$ with the interpretation that $\delta(\mathbf{a}|\mathbf{z})$ is the conditional probability of action $\mathbf{a}$ given $\mathbf{z}$ is observed. For each $\mathbf{z}$, a nonrandomized decision rule chooses a single element of $\Gamma$ with probability 1 and assigns all other actions probability 0.

Let $\psi_i(\mathbf{z})$ be the probability of rejecting $H_i$. A decision procedure $\delta(\mathbf{a}|\mathbf{z})$ determines a set of $\psi_i^{(\delta)}(\mathbf{z})$, $i=1,\dots,k$, as follows:

$$\psi_i^{\delta}(\mathbf{z})=\sum_{\mathbf{a}\in A_i}\delta(\mathbf{a}|\mathbf{z}), \tag{2.1}$$

where $A_i=\{\mathbf{a}\in\Gamma:\mathbf{a}\text{ has a 1 in the }i\text{th position}\}$. Whereas $\delta(\mathbf{a}|\mathbf{z})$ determines $\boldsymbol{\psi}(\mathbf{z})$, the reverse is not true. If $\boldsymbol{\psi}(\mathbf{z})=(\psi_1,\dots,\psi_k)'$ is nonrandomized, it uniquely determines some $\delta(\mathbf{a}|\mathbf{z})$. The $\delta(\mathbf{a}|\mathbf{z})$ determined is also nonrandomized.

The parameter space is $\Omega=\{\boldsymbol{\mu}:\mu_i\ge 0,i=1,\dots,k\}$. Partition the parameter space $\Omega$ into $2^k$ sets $\Omega_{\mathbf{v}}$, $\mathbf{v}\in\Gamma$, where $\Omega_{\mathbf{v}}=\{\boldsymbol{\mu}:\mu=(\mu_1,\mu_2,\dots,\mu_k)'$, $\mu_i>0$ if $v_i=1$ and $\mu_i=0$ if $v_i=0$, $i=1,\dots,k\}$. Also let $\Omega^{(i)}=\{\boldsymbol{\mu}:\boldsymbol{\mu}\in\Omega,\mu_i=0\}$ and let $\Omega^{(i)c}$ be the complement of $\Omega^{(i)}$ relative to $\Omega$.

A loss function is a function of the action taken and the true state of nature. We will study several different loss functions and their relationships. For each individual hypothesis $H_i$ the loss function is zero for a correct decision, 1 for rejecting $H_i$ when it is true and 1 for accepting $H_i$ when it is false. (Note that the ensuing development would also work if 1 is replaced by $b$, $b>0$, when a false acceptance is made.) Such a loss function determines a risk

$$R_{(i)}(\psi_i,\boldsymbol{\mu})=(1-v_i)E_{\boldsymbol{\mu}}\psi_i(z)+v_i(1-E_{\boldsymbol{\mu}}\psi_i(z)). \tag{2.2}$$

For the finite action problem a sum of the loss functions of the individual problems is

$$L(\mathbf{a},\boldsymbol{\mu})=\sum_{i=1}^k a_i(1-v_i)+\sum_{i=1}^k(1-a_i)v_i,\qquad \boldsymbol{\mu}\in\Omega_{\mathbf{v}}. \tag{2.3}$$

The corresponding risk function is $\sum_{i=1}^{k} R_{(i)}$, which can be expressed as

$$E_{\boldsymbol{\mu}}(\boldsymbol{\psi}'(1-\mathbf{v}) + (1-\boldsymbol{\psi})'\mathbf{v}). \tag{2.4}$$

The risk function (2.4) is used by Lehmann (1957) and CS (2005). CS (2005) also study a vector loss function consisting of the vector of losses for the individual component problems. The corresponding vector risk called VRI is

$$(R_{(1)}(\boldsymbol{\psi}_1, \boldsymbol{\mu}), \ldots, R_{(k)}(\boldsymbol{\psi}_k, \boldsymbol{\mu})). \tag{2.5}$$

Another vector loss consisting of two components $(L_0, L_1)$ is where

$$L_0(\mathbf{a}, \boldsymbol{\mu}) = \sum_{i=1}^{k} a_i(1-v_i), \qquad L_1(\mathbf{a}, \boldsymbol{\mu}) = \sum_{i=1}^{k} (1-a_i)v_i, \qquad \boldsymbol{\mu} \in \Omega_{\mathbf{v}}. \tag{2.6}$$

The corresponding risk function can be expressed as $(R_0(\boldsymbol{\psi}, \boldsymbol{\mu}), R_1(\boldsymbol{\psi}, \boldsymbol{\mu}))$, where

$$R_0(\boldsymbol{\psi}, \boldsymbol{\mu}) = \sum_{i=1}^{k} (1-v_i) E_{\boldsymbol{\mu}} \psi_i(\mathbf{z}), \qquad R_1(\boldsymbol{\psi}, \boldsymbol{\mu}) = \sum_{i=1}^{k} v_i E_{\boldsymbol{\mu}}(1-\psi_i(\mathbf{z})). \tag{2.7}$$

Suppose $m$ is the number of positive $\mu_i$. Then according to the definition of average power given by Benjamini and Hochberg (1995), and noted by Shaffer (1995) and Dudoit, Shaffer and Boldrick (2003), $R_1/m$ is 1 minus the average power. $R_0/(k-m)$ is the average size. Also one may justifiably call $R_0$ or $R_0/(k-m)$ the false detection rate or the false rejection rate (FRR) and call $R_1$ or $R_1/m$ the false acceptance rate (FAR). We call the vector risk $(R_0, R_1)$ in (2.7) VRSP since it is related to average size and average power. We note that the class of admissible procedures is largest for the VRI risk function in (2.5) and smallest for the risk function in (2.4). Yet the class of admissible procedures is certainly larger for VRSP than for the risk function in (2.4). Thus any procedure which is inadmissible for (2.7) is also inadmissible for (2.4).

In this paper we focus on VRSP. To deal with VRSP we use a device utilized by Cohen and Sackrowitz (1984). That is, we introduce an artificial but useful problem. Let $\theta$ be a nuisance parameter which takes on the value 0 or 1. Next define the one-dimensional loss function

$$L^*(\mathbf{a}, (\boldsymbol{\mu}, \theta)) = L_\theta(\mathbf{a}, \boldsymbol{\mu}). \tag{2.8}$$

It now follows from Cohen and Sackrowitz (1984) that the class of admissible procedures for the problem using (2.6) as a loss function is the same as the problem using (2.8) as a loss function but treating $\theta$ as a parameter which can either be 0 or 1. Hence we study the problem using (2.8) as the loss function. The corresponding risk function will be denoted as $R^*(\boldsymbol{\psi}, (\boldsymbol{\mu}, \theta))$.

Now a decision procedure $\boldsymbol{\psi}^*$ is Bayes with respect to (w.r.t.) a prior distribution $\xi(\boldsymbol{\mu}, \theta)$ if

$$E_\xi R^*(\boldsymbol{\psi}^*, (\boldsymbol{\mu}, \theta)) = \inf_{\boldsymbol{\psi}} E_\xi R^*(\boldsymbol{\psi}, (\boldsymbol{\mu}, \theta)). \tag{2.9}$$

The prior distribution is written as

$$\xi(\boldsymbol{\mu}, \theta) = \begin{cases} \xi_0(\boldsymbol{\mu})\beta, & \text{if } \theta = 0, \\ \xi_1(\boldsymbol{\mu})(1-\beta), & \text{if } \theta = 1, \end{cases}$$

where $\beta$ is the probability that $\theta = 0$ and $\xi_0(\boldsymbol{\mu})$ is the conditional distribution of $\boldsymbol{\mu}$ given $\theta = 0$ and where $(1-\beta)$ is the probability that $\theta = 1$ and $\xi_1(\boldsymbol{\mu})$ is the conditional distribution of $\boldsymbol{\mu}$ given $\theta = 1$. We write the density of $\mathbf{z}$ given $(\boldsymbol{\mu}, \theta)$, noting that this density is the same regardless of $\theta$. That is, $f(\mathbf{z}|\boldsymbol{\mu}, 0) = f(\mathbf{z}|\boldsymbol{\mu}, 1) = f(\mathbf{z}|\boldsymbol{\mu})$ where

$$f(\mathbf{z}|\boldsymbol{\mu}) = (1/(2\pi)|\Sigma|^{1/2}) e^{-(1/2)z'\Sigma^{-1}z} e^{z'\Sigma^{-1}\boldsymbol{\mu}} e^{-(1/2)\boldsymbol{\mu}'\Sigma^{-1}\boldsymbol{\mu}}. \tag{2.10}$$

Note the marginal distribution of $\mathbf{z}$ is

$$f(\mathbf{z}) = \int_\Omega f(\mathbf{z}|\boldsymbol{\mu})[\beta\, d\xi_0(\boldsymbol{\mu}) + (1-\beta)\, d\xi_1(\boldsymbol{\mu})]. \tag{2.11}$$

The following theorem describes a Bayes procedure.

Theorem 2.1. *Consider the risk function $R^*(\boldsymbol{\psi}, (\boldsymbol{\mu}, \theta))$. The Bayes procedure w.r.t. $\xi(\boldsymbol{\mu}, \theta)$ is $\boldsymbol{\psi}^* = (\psi_1^*, \ldots, \psi_i^*)'$ where*

$$\psi_i^* = \begin{cases} 1, & \text{if } \dfrac{\int_{\Omega^{(i)}} f(\mathbf{z}|\boldsymbol{\mu})[\beta\, d\xi_0(\boldsymbol{\mu}) + (1-\beta)\, d\xi_1(\boldsymbol{\mu})]}{\int_\Omega f(\mathbf{z}|\boldsymbol{\mu})\, d\xi_1(\boldsymbol{\mu})} < (1-\beta), \\ 0, & \text{otherwise.} \end{cases} \tag{2.12}$$

Proof. The risk function $R^*(\boldsymbol{\psi}, (\boldsymbol{\mu}, \theta))$ can be written as

$$\begin{aligned} R^*(\boldsymbol{\psi}, (\boldsymbol{\mu}, \theta)) &= (1-\theta)\sum_{i=1}^k (1-v_i) E_{\boldsymbol{\mu}}[\psi_i(\mathbf{z})] + \theta \sum_{i=1}^k v_i E_{\boldsymbol{\mu}}[1-\psi_i(\mathbf{z})] \\ &= \sum_{i=1}^k \{(\theta v_i) + (1-\theta-v_i) E_{\boldsymbol{\mu}}\psi_i(\mathbf{z})\}, \qquad \boldsymbol{\mu} \in \Omega_{\mathbf{v}}. \end{aligned} \tag{2.13}$$

To find the Bayes procedure we must minimize the expected risk. Using (2.13) we see that this amounts to setting $\psi_i(\mathbf{z}) = 1$ if the posterior expected value

$$E\{1 - \Theta - V_i|\mathbf{z}\} < 0, \tag{2.14}$$

where now $\Theta$ and $V_i$ are regarded as random variables with joint prior distribution determined by $\xi(\boldsymbol{\mu}, \theta)$. The left-hand side of (2.14) reduces to

$$P\{V_i = 0|\mathbf{z}\} - P\{\Theta = 1|\mathbf{z}\}. \tag{2.15}$$

Now (2.15) is

$$\begin{aligned}\Bigg\{\int_{\Omega^{(i)}} f(\mathbf{z}|\boldsymbol{\mu})[\beta\, d\,\xi_0(\boldsymbol{\mu}) + (1-\beta)\, d\xi_1(\boldsymbol{\mu})]\\ - (1-\beta)\int_{\Omega} f(\mathbf{z}|\mu)\, d\xi_1(\boldsymbol{\mu})\Bigg\}\Big/ f(\mathbf{z}).\end{aligned} \tag{2.16}$$

We see that (2.14) and (2.16) lead to (2.12). □

The step-up procedure we study is as follows:

PROCEDURE 2.1. Let $Z_{(1)} \le Z_{(2)} \le \cdots \le Z_{(k)}$ be the order statistics for the set $(Z_1, \dots, Z_k)$ and let $C_j$ be a strictly increasing set of critical values:

(i) If $Z_{(1)} \le C_1$, accept $H_{(1)}$ where $H_{(1)}$ is the hypothesis corresponding to $Z_{(1)}$. Otherwise reject all $H_i$.

(ii) If $H_{(1)}$ is accepted, accept $H_{(2)}$ if $Z_{(2)} \le C_2$. Otherwise reject $H_{(2)}, \dots, H_{(k)}$.

(iii) In general, at stage $j$, if $Z_{(j)} \le C_j$, accept $H_{(j)}$. Otherwise reject $H_{(j)}, \dots, H_{(k)}$.

Call the step-up procedure $\boldsymbol{\psi}_{\mathrm{SU}}(\mathbf{z})$. The procedure for $k = 2$ is shown in Figure 1.

**3. Characterization of symmetric Bayes procedures.** In order to characterize symmetric Bayes procedures we first recognize that the problem with loss function (2.8) is invariant under the following groups of transformations:

(i) $G = \{g : g\mathbf{z}$ is a permutation of the coordinates of $\mathbf{z}$; i.e., $g$ is a $k \times k$ permutation matrix$\}$.

(ii) $\overline{G} = \{\bar{g} : \bar{g}(\boldsymbol{\mu}, \theta)$ is a permutation of the coordinates of $\boldsymbol{\mu}$ while leaving $\theta$ as is; i.e., $\bar{g} = \begin{pmatrix} g & \mathbf{0} \\ \mathbf{0} & 1 \end{pmatrix}\}$.

(iii) $\widetilde{G} = \{\tilde{g} : \tilde{g}(\mathbf{a})$ is a permutation of the coordinates of $\mathbf{a}$, i.e., $\tilde{g} = g\}$.

Since the problem is invariant under the finite group $G$, it follows from Ferguson [(1967), Theorem 3, page 162] that any symmetric Bayes procedure is Bayes w.r.t. an invariant prior distribution. Any invariant prior distribution (under $\bar{g}$) depends only on the maximal invariant parameter $(\mu_{(1)}, \dots, \mu_{(k)}, \theta)$. This restriction then implies that for Bayes procedures, all prior distributions are symmetric in $\boldsymbol{\mu}$ for each fixed $\theta$. In particular, the conditional distributions $\xi_0(\boldsymbol{\mu})$ and $\xi_1(\boldsymbol{\mu})$ will be permutation invariant.

In order to characterize symmetric Bayes procedures we will be using (2.12). We first wish to express the integrand in a simplified fashion. Toward this end recall the expression (2.10) for $f(\mathbf{z}|\boldsymbol{\mu})$. Since $\Sigma$ is intraclass, that is, $\Sigma = \sigma^2(1-\rho)I + \rho\mathbf{1}\mathbf{1}'$, $\mathbf{1} = (1,\dots,1)'$, $\Sigma^{-1} = (\sigma^2(1-\rho))^{-1}(I - G11')$ where $G = \rho/(1+(k-1)\rho)$, we can express the numerator of (2.12) as

$$e^{-(1/2)\mathbf{z}'\Sigma^{-1}\mathbf{z}} \int_{\Omega^{(i)}} \exp(\mathbf{z}'\Sigma^{-1}\boldsymbol{\mu} - (1/2)\boldsymbol{\mu}'\Sigma^{-1}\boldsymbol{\mu}) \times [\beta\, d\xi_0(\boldsymbol{\mu}) + (1-\beta)\, d\xi_1(\boldsymbol{\mu})]. \tag{3.1}$$

Noting that $\mathbf{z}'\Sigma^{-1}\boldsymbol{\mu} = (\sigma^2(1-\rho))^{-1}(\mathbf{z}'\boldsymbol{\mu} - G\mathbf{1}'\mathbf{z}\mathbf{1}'\boldsymbol{\mu})$, letting for fixed $\mathbf{z}$,

$$d\xi_\theta^*(\boldsymbol{\mu}) = \exp\{-(\sigma^2(1-\rho))^{-1}G\mathbf{1}'\mathbf{z}\mathbf{1}'\boldsymbol{\mu} + (1/2)\boldsymbol{\mu}'\Sigma^{-1}\boldsymbol{\mu}\}\, d\xi_\theta(\boldsymbol{\mu}), \tag{3.2}$$

and without loss of generality taking $\sigma^2(1-\rho) = 1$, we can rewrite (2.12) as

$$\psi_i^* = \begin{cases} 1, & \text{if } Q(\Omega^{(i)}|\mathbf{z}) = \left\{\displaystyle\int_{\Omega^{(i)}} e^{\mathbf{z}'\boldsymbol{\mu}}[\beta\, d\xi_0^*(\boldsymbol{\mu}) + (1-\beta)\, d\xi_1^*(\boldsymbol{\mu})] \times \left(\int_\Omega e^{\mathbf{z}'\boldsymbol{\mu}}\, d\xi_1^*(\boldsymbol{\mu})\right)^{-1}\right\} \\ & \qquad < 1-\beta, \\ 0, & \text{otherwise.} \end{cases} \tag{3.3}$$

Note that in (3.2) we absorb an expression involving $\mathbf{z}$ into the prior. This is okay since $\mathbf{z}$ is fixed and in the development to follow even when $\mathbf{z}$ changes $\mathbf{1}'\mathbf{z}$ will remain constant.

To characterize symmetric Bayes rules it suffices to consider only sample points such that $z_1 \le \cdots \le z_k$. Now we give:

Lemma 3.1. *Fix* $\mathbf{z}$ *and assume* $z_1 \le \cdots \le z_k$. *Then*

$$Q(\Omega^{(i)}|\mathbf{z}) \ge Q(\Omega^{(i+1)}|\mathbf{z}), \qquad i = 1,\dots,k-1. \tag{3.4}$$

*The inequality is strict unless* $z_i = z_{i+1}$.

Proof. We need only consider the integral in the numerator of (3.3). Writing $[\beta\, d\xi_0^*(\boldsymbol{\mu}) + (1-\beta)\, d\xi_1^*(\boldsymbol{\mu})]$ as $d\xi^*(\boldsymbol{\mu})$, we note

$$\int_{\Omega^{(i)}} e^{\mathbf{z}'\boldsymbol{\mu}}\, d\xi^*(\boldsymbol{\mu}) = \int_{\Omega^{(i)}} \exp\left(\sum_{j\ne i,i+1} z_j\mu_j + z_{i+1}\mu_{i+1}\right) d\xi^*(\boldsymbol{\mu}). \tag{3.5}$$

Make the change of variables $\mu_i = \mu_{i+1}$, $\mu_{i+1} = \mu_i$ in (3.5) to find (3.5) is equal to

$$\begin{aligned} &\int_{\Omega^{(i+1)}} \exp\left(\sum_{j\ne i,i+1} z_j\mu_j + z_{i+1}\mu_i\right) d\xi^*(\boldsymbol{\mu}) \\ &\qquad \ge \int_{\Omega^{(i+1)}} \exp\left(\sum_{j\ne i,i+1} z_j\mu_j + z_i\mu_i\right) d\xi^*(\boldsymbol{\mu}). \end{aligned} \tag{3.6}$$

Thus from (3.5) and (3.6) we have (3.4). Note that the inequality in the proof is strict unless $z_i = z_{i+1}$. This completes the proof of the lemma. □

THEOREM 3.2. *Let* $\mathbf{z}$ *be such that* $z_1 < \cdots < z_k$. *Let* $r \in \{0, 1, \ldots, k\}$ *be the element of the set for which* $Q(\Omega^{(r)}|\mathbf{z}) > (1-\beta) > Q(\Omega^{(r+1)}|\mathbf{z})$, *where* $r = 0$ *means* $Q(\Omega^{(i)}|\mathbf{z}) < (1-\beta)$ *for all* $i = 1, \ldots, k$ *and* $r = k$ *means* $Q(\Omega^{(i)}|\mathbf{z}) > (1-\beta)$ *for all* $i = 1, \ldots, k$. *Then the Bayes procedure is* $\psi_i(\mathbf{z}) = 0$, $i = 1, \ldots, r$, $\psi_i(\mathbf{z}) = 1$, $i = r+1, \ldots, k$.

PROOF. Use Theorem 2.1 and Lemma 3.1. □

**4. Complete class and inadmissibility of step-up.** Symmetric Bayes procedures and weak * limits of sequences of symmetric Bayes procedures against symmetric prior distributions form a complete class of symmetric procedures for this problem. See Weiss [(1961), page 81], where he defines a weak * limit as follows: Let $\boldsymbol{\psi}_n$ be a sequence of procedures. Then $\boldsymbol{\psi}_n$ converges to $\boldsymbol{\psi}$ if

$$\lim_{n\to\infty} R(\boldsymbol{\psi}_n, \boldsymbol{\mu}) = R(\boldsymbol{\psi}, \boldsymbol{\mu}).$$

Another complete class of procedures for this problem is the set of almost everywhere (a.e.) nonrandomized procedures. This follows from a result in Matthes and Truax (1967) where it is demonstrated that each admissible $\psi_i(\mathbf{z})$ must be nonrandomized a.e. It follows that the nonrandomized symmetric Bayes procedures and their a.e. limits are a complete class of symmetric procedures for this problem.

We proceed to give a necessary condition for admissibility based on a complete class. Toward this end let $t_j$ be the following partial sums of $(z_1, \ldots, z_k)$. That is, let $t_j = \sum_{i=j}^{k} z_i$, $j = 1, \ldots, k$. Let $t_{k+1} = 0$ and $t_0 = -\infty$.

LEMMA 4.1. *Let* $\mathcal{S} = \{\mathbf{t} : t_k > t_{k-1} - t_k > \cdots > t_1 - t_2\}$. *Then for* $j = 2, \ldots, k$, $\mathbf{t} \in \mathcal{S}$, $Q(\Omega^{(j)}|\mathbf{t})$ *as a function of* $t_j$ *is strictly decreasing while* $t_1, \ldots, t_{j-1}, t_{j+1}, \ldots, t_k$ *are held fixed.*

PROOF. Note we may write

$$Q(\Omega^{(j)}|\mathbf{t}) = \frac{\int_{\Omega^{(j)}} \exp(\sum_{i=1, i\neq j}^{k} (t_i - t_{i+1})\mu_i)\, d\xi^*(\boldsymbol{\mu})}{\int_{\Omega} \exp(\sum_{i=1}^{k} (t_i - t_{i+1})\mu_i)\, d\xi_1^*(\boldsymbol{\mu})}. \tag{4.1}$$

For fixed $t_1, t_2, \ldots, t_{j-1}, t_{j+1}, \ldots, t_k$ the numerator is a strictly decreasing function of $t_j$ (recall $\mu_{j-1} \geq 0$) while the denominator, being a Schur convex function of $\mathbf{z}$ (it is convex and permutation invariant), is an increasing function of $t_j$, $t_j \in \mathcal{S}$, while all other partial sums are fixed. See Marshall and Olkin (1979) for discussion of Schur convex functions. It follows that $Q(\Omega^{(j)}|\mathbf{t})$ then is a decreasing function of $t_j$. This completes the proof of the lemma. □

LEMMA 4.2. *Let $j = 2, \ldots, k$. Let $\boldsymbol{\psi}(\mathbf{t})$ be a symmetric Bayes procedure. Then for $\mathbf{t} \in \mathcal{S}$, $\psi_j(\mathbf{t})$ is a nondecreasing function of $t_j$ while $(t_1, \ldots, t_{j-1}, t_{j+1}, \ldots, t_k)$ are fixed.*

PROOF. Note since $\boldsymbol{\psi}(\mathbf{t})$ is a symmetric Bayes procedure it follows from the proof of Lemma 4.1 that $\psi_j(\mathbf{t})$ is nonrandomized for $j = 2, \ldots, k$. Use Lemma 4.1 again to conclude that for $\mathbf{t} \in \mathcal{S}$, $\psi_j(\mathbf{t})$ is a nondecreasing function of $t_j$ while $(t_1, \ldots, t_{j-1}, t_{j+1}, \ldots, t_k)$ are fixed. □

THEOREM 4.3. *Let $j = 1, \ldots, k-1$. Let $\boldsymbol{\psi}(\mathbf{t})$ be a symmetric procedure such that there exists a sample point $\mathbf{t}^* \in \mathcal{S}$ for which $\psi_j(\mathbf{t}^*) = 0$. Then a necessary condition for $\boldsymbol{\psi}(\mathbf{t})$ to be admissible is that $\psi_j(\mathbf{t}) = 0$ for all $\mathbf{t} \in \mathcal{S}$ such that $t_j < t_j^*$.*

PROOF. Recall that symmetric Bayes and a.e. limits of sequences of symmetric Bayes procedures are a complete class of symmetric procedures. Now Lemma 4.2 implies that every Bayes procedure has the required property. The required property must also hold for any a.e. limit of a sequence of symmetric Bayes procedures. To see this let $\boldsymbol{\psi}_n(\mathbf{t}) = (\psi_{1n}(t), \ldots, \psi_{kn}(t))'$ be a sequence of symmetric Bayes procedures with $\boldsymbol{\psi}(\mathbf{t})$ its a.e. limit. Since $\psi_{jn}(\mathbf{t})$ is a nondecreasing function it follows that its a.e. limit is also a nondecreasing function. This establishes the theorem. □

COROLLARY 4.4. *Let $\boldsymbol{\psi}(\mathbf{z})$ be a procedure such that there exists a sample point $\mathbf{z}^* = (z_1^*, \ldots, z_k^*)'$ with $z_k^* > z_{k-1}^* > \cdots > z_1^*$ for which $\psi_k(\mathbf{z}^*) = 0$. Then a necessary condition for $\boldsymbol{\psi}(\mathbf{z})$ to be admissible is that $\psi_k(\mathbf{z}) = 0$ for all $\mathbf{z}$ in the set $\{\mathbf{z} : (z_k^* + z_{k-1}^*)/2 \le z_k \le z_k^*, z_{k-2} = z_{k-2}^*, \ldots, z_1 = z_1^*, z_k + z_{k-1} = z_k^* + z_{k-1}^*\}$.*

PROOF. This follows from Theorem 4.3 since $t_k = z_k$, and fixing $t_{k-1}, z_{k-2}, \ldots, z_1$ is equivalent to fixing $t_{k-1}, \ldots, t_1$. □

THEOREM 4.5. *The step-up procedure given in Procedure* 2.1 *is inadmissible.*

PROOF. We show that Procedure 2.1 does not satisfy the necessary condition for admissibility given in Corollary 4.4. Consider the sample point $\mathbf{z}^* = (z_1^*, z_2^*, \ldots, z_k^*)'$ where $z_j^* = C_j - \varepsilon$, $j = 1, \ldots, k$, for some $\varepsilon > 0$ to be chosen. Note since $C_1 < \cdots < C_k$, $z_1^* < z_2^* < \cdots < z_k^*$ and also note that $\boldsymbol{\psi}_{\mathrm{SU}}(\mathbf{z}^*) = \mathbf{0}$. In particular, the last coordinate of $\boldsymbol{\psi}_{\mathrm{SU}}(\mathbf{z}^*)$ is zero. Now consider the sample point $\bar{\mathbf{z}} = (z_1^*, z_2^*, \ldots, z_{k-2}^*, \bar{z}_{k-1}, \bar{z}_k)$ where $\bar{z}_k = \bar{z}_{k-1} = [(C_k + C_{k-1})/2] - \varepsilon$. Notice that for sufficiently small $\varepsilon$, $\bar{z}_{k-1} > C_{k-1}$, which means that $\boldsymbol{\psi}_{\mathrm{SU}}(\bar{z}) = (0, 0, \ldots, 0, 1, 1)'$. In fact there is an open interval of $\mathbf{z}$ points on the line $z_k + z_{k-1} = z_k^* + z_{k-1}^* = t_{k-1}^*, z_{k-2}^*, \ldots, z_1^*$ beginning at $z_k = t_{k-1}^*/2$ and ending before $z_k = z_k^*$ such that $\boldsymbol{\psi}_{\mathrm{SU}}(\mathbf{z}) = (0, 0, \ldots, 0, 1, 1)'$. In particular, the last coordinate of $\boldsymbol{\psi}_{\mathrm{SU}}(\mathbf{z}) = 1$. This represents a violation of the necessary condition for admissibility given in Corollary 4.4.

The result of Theorem 4.5 is, in a sense, akin to the famous inadmissibility result of Stein (1956). Stein considered the model $\mathbf{Z} \sim N(\boldsymbol{\mu}, I)$ and proved that if the loss function is the sum of squared errors, then $\mathbf{Z}$ is an inadmissible estimator of $\boldsymbol{\mu}$ when $k \geq 3$. This in spite of the fact that each $Z_i$ is admissible for $\mu_i$ if the loss function is squared error. In our multiple endpoints testing problem $\boldsymbol{\psi}_{\mathrm{SU}} = (\psi_{\mathrm{SU}1}, \psi_{\mathrm{SU}2}, \ldots, \psi_{\mathrm{SU}k})'$ is such that $\psi_{\mathrm{SU}i}$ is an admissible test of $H_i : \mu_i = 0$ vs $K_i : \mu_i > 0$ when the loss function is $(0, 1)$ and $\rho \geq 0$. See CS (2005). Yet $\boldsymbol{\psi}_{\mathrm{SU}}$ is inadmissible as a finite action procedure when the loss is the sum of losses of the component problems [or for the vector loss (2.6)]. Here the result is true for $k \geq 2$. $\square$

**5. A procedure which beats step-up.** In the case of $k = 2$, the step-up procedure is shown in Figure 1(a). It is easily seen that the necessary condition of Theorem 4.3 is violated when $2C_1 < Z_1 + Z_2 < C_1 + C_2$. This is the shaded strip in Figure 1(b). By making changes in this strip we show how to construct a procedure, $\boldsymbol{\psi}^*$, that has a vector risk which is less than or equal to the risk of step-up for all $\boldsymbol{\mu}$.

We begin with any fixed $t \in (2C_1, C_1 + C_2)$ and consider $\mathbf{Z}$ such that $Z_1 + Z_2 = t$. Without loss of generality let $\sigma^2 = 1$. We note that the conditional distribution of $Z_1$ given $Z_1 + Z_2 = t$ is $N(\frac{t}{2} + \frac{1}{2}(\mu_1 - \mu_2), \frac{1}{2}(1 - \rho))$. Also (see Figure 2), for the step-up procedure $\boldsymbol{\psi}_{\mathrm{SU}}(\mathbf{z})$, when $z_1 + z_2 = t$, we have

$$\boldsymbol{\psi}_{\mathrm{SU}}(\mathbf{z}) = \begin{cases} (0, 1), & \text{if } z_1 < t - C_2, \\ (0, 0), & \text{if } t - C_2 < z_1 < C_1, \\ (1, 1), & \text{if } C_1 < z_1 < t - C_1, \\ (0, 0), & \text{if } t - C_1 < z_1 < C_2, \\ (1, 0), & \text{if } C_2 < z_1. \end{cases}$$

The procedure $\boldsymbol{\psi}^*$ is constructed as follows. Consider $P_{\mu_1 = \mu_2}(t - C_1 < Z_1 < C_2 | Z_1 + Z_2 = t) - P_{\mu_1 = \mu_2}(\frac{t}{2} < Z_1 < t - C_1 | Z_1 + Z_2 = t)$. If we let $D(t)$ be this difference in conditional probabilities, then

$$D(t) = \Phi\left(\frac{2C_2 - t}{\sqrt{2(1 - \rho)}}\right) - 2\Phi\left(\frac{t - 2C_1}{\sqrt{2(1 - \rho)}}\right) + \frac{1}{2}. \tag{5.1}$$

Next define $C^* = C^*(t)$ by setting

$$P\left(\frac{t}{2} < Z_1 < C^* | Z_1 + Z_2 = t\right) = \Phi\left(\frac{2C^* - t}{\sqrt{2(1 - \rho)}}\right) - \frac{1}{2} = |D(t)|. \tag{5.2}$$

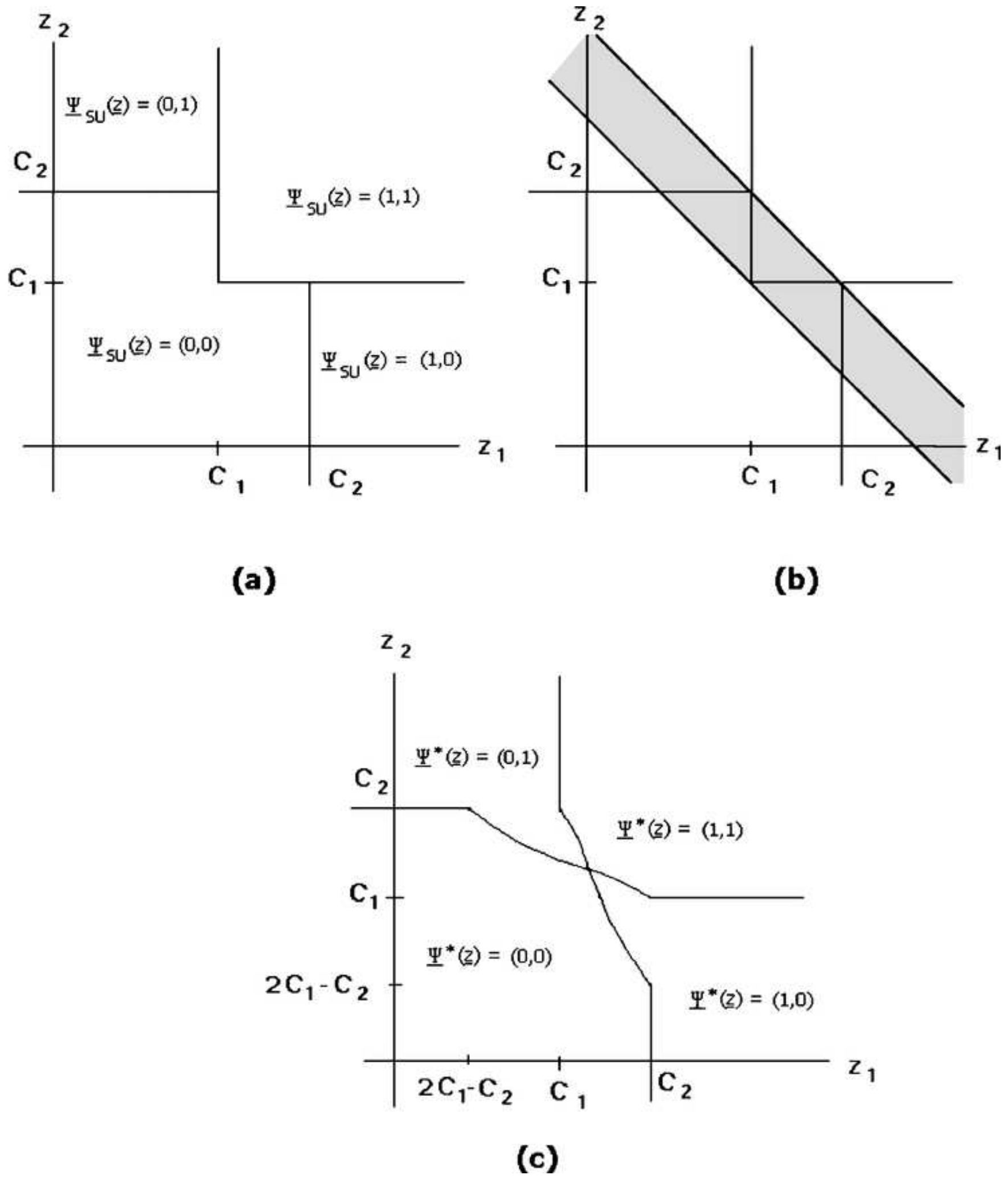

FIG. 1. (a) *and* (b) *show the step-up procedure* $\phi_{\mathrm{SU}}$, (c) *shows* $\phi$.

That is, $C^*$ is the solution to

$$\Phi\left(\frac{2C^* - t}{\sqrt{2(1-\rho)}}\right) - \frac{1}{2} = \left|\frac{1}{2} + \Phi\left(\frac{2C_2 - t}{\sqrt{2(1-\rho)}}\right) - 2\Phi\left(\frac{t - 2C_1}{\sqrt{2(1-\rho)}}\right)\right|. \tag{5.3}$$

Then for $z_1 + z_2 = t$ and $D(t) > 0$, let

$$\boldsymbol{\psi}^*(z_1, z_2) = \begin{cases} (0,1), & \text{if } z_1 < t - C^*, \\ (0,0), & \text{if } t - C^* < z_1 < C^*, \\ (1,0), & \text{if } C^* < z_1. \end{cases}$$

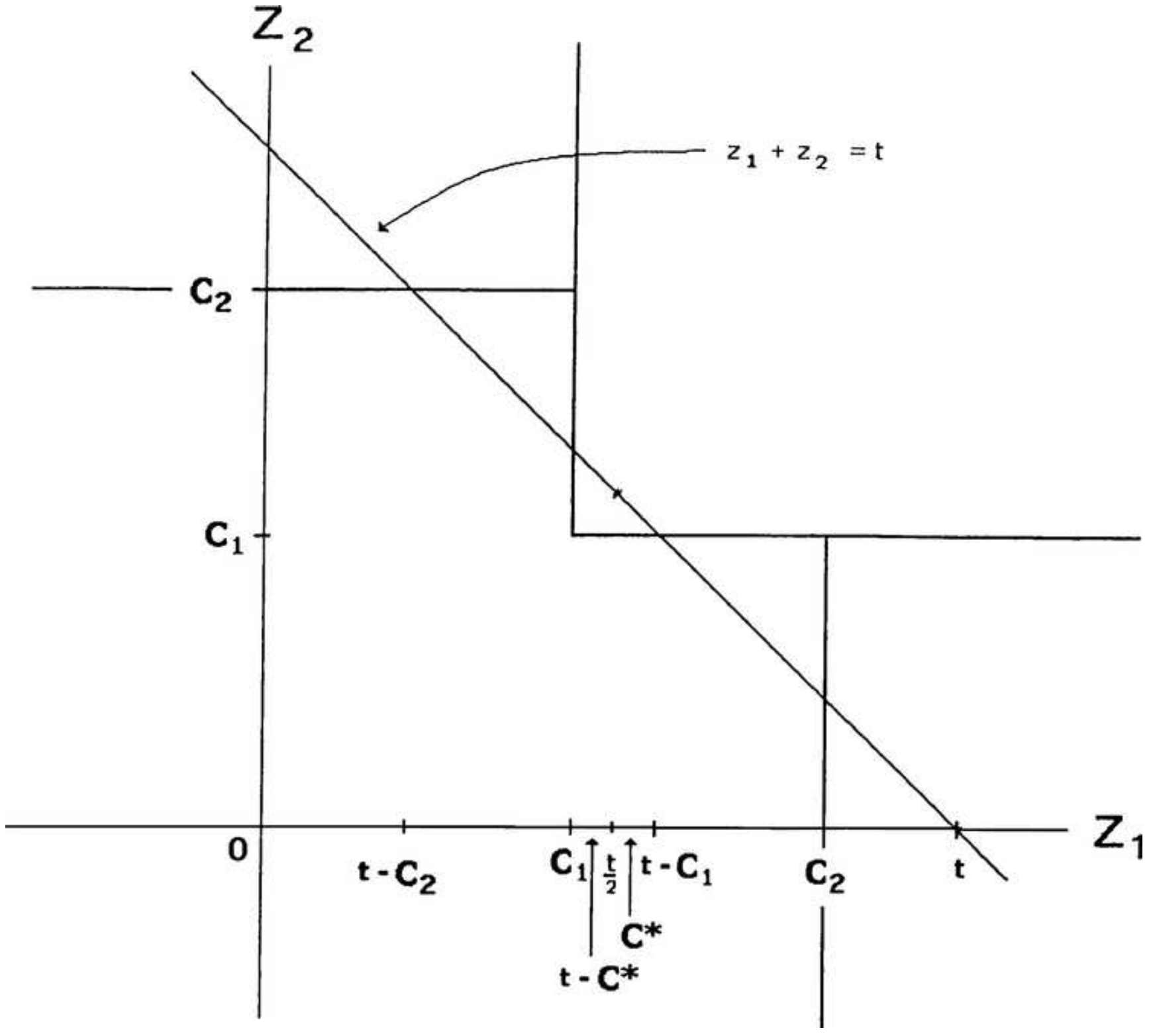

FIG. 2. *The step-up procedure with the line $z_1 + z_2 = t$.*

On the other hand, if $D(t) < 0$, let

$$\boldsymbol{\psi}^*(z_1, z_2) = \begin{cases} (0,1), & \text{if } z_1 < t - C^*, \\ (1,1), & \text{if } t - C^* < z_1 < C^*, \\ (1,0), & \text{if } C^* < z_1. \end{cases}$$

The resulting procedure is sketched in Figure 1(c).

THEOREM 5.1. *The procedure $\boldsymbol{\psi}^*$ is better than $\boldsymbol{\psi}_{\text{SU}}$ for the vector risk VRSP.*

PROOF. If we let $\boldsymbol{\psi}_{\text{SU}}(\mathbf{z})$ denote the step-up procedure, then it will be shown that the procedure $\boldsymbol{\psi}^*(\mathbf{z})$ above is such that

$$R_0(\boldsymbol{\psi}^*, \boldsymbol{\mu}) + bR_1(\boldsymbol{\psi}^*, \boldsymbol{\mu}) \le R_0(\boldsymbol{\psi}_{\text{SU}}, \boldsymbol{\mu}) + bR_1(\boldsymbol{\psi}_{\text{SU}}, \boldsymbol{\mu}), \tag{5.4}$$

with strict inequality for some $\boldsymbol{\mu}$, for every $b > 0$. Note that the procedure $\boldsymbol{\psi}^*$ does not depend on $b$. This implies that $\boldsymbol{\psi}^*$ beats $\boldsymbol{\psi}_{\text{SU}}$ for VRSP.

Using (2.7) we show (5.4) by showing

$$E_{\boldsymbol{\mu}}\{(\boldsymbol{\psi}_{\text{SU}}(\mathbf{z}) - \boldsymbol{\psi}^*(\mathbf{z}))'(\mathbf{1} - (b+1)\mathbf{v}) | Z_1 + Z_2 = t\} > 0 \tag{5.5}$$

TABLE 1

*Evaluation of $\phi_{\mathrm{SU}}(z)$, $\phi^*(\mathbf{z})$ and $W(\mathbf{z};\mathbf{v})$ when $z_1+z_2=t$*

| | | | | $\mathbf{v}$ | | | |
|---|---|---|---|---|---|---|---|
| | | | | $(0,0)$ | $(1,0)$ | $(0,1)$ | $(1,1)$ |
| | | | | $\mathbf{1}-(b+1)\mathbf{v}$ | | | |
| | | | | $(1,1)$ | $(-b,1)$ | $(1,-b)$ | $(-b,-b)$ |
| | | | | $W(\mathbf{z};\mathbf{v})$ | | | |
| | $\phi_{\mathrm{SU}}(\mathbf{z})'$ | $\phi^*(\mathbf{z})'$ | $(\phi_{\mathrm{SU}}(\mathbf{z})'-\phi^*(\mathbf{z}))'$ | $(\phi_{\mathrm{SU}}(\mathbf{z})-\phi^*(\mathbf{z}))'(\mathbf{1}-(b+1)\mathbf{v})$ | | | |
| $-\infty\le z_1\le t-C_2$ | (0, 1) | (0, 1) | (0,0) | 0 | 0 | 0 | 0 |
| $t-C_2\le z_1\le C_1$ | (0, 0) | (0,1) | (0, −1) | −1 | −1 | $b$ | $b$ |
| $C_1\le z_1\le t-C^*$ | (1, 1) | (0,1) | (1,0) | 1 | $-b$ | 1 | $-b$ |
| $t-C^*\le z_1\le C^*$ | (1, 1) | (0, 0) | (1,1) | 2 | $1-b$ | $1-b$ | $-2b$ |
| $C^*\le z_1\le t-C_1$ | (1, 1) | (1, 0) | (0,1) | 1 | 1 | $-b$ | $-b$ |
| $t-C_1\le z_1\le C_2$ | (0, 0) | (1, 0) | ( − 1,0) | −1 | $b$ | −1 | $b$ |
| $C_2\le z_1\le\infty$ | (1, 0) | (1, 0) | (0,0) | 0 | 0 | 0 | 0 |

for all $\boldsymbol{\mu}\in\Omega_{\mathbf{v}}$, $\mathbf{v}\in\Gamma$, and $t\in(2C_1, C_1+C_2)$. We will only study the case of $D(t)>0$ and $\frac{1}{2}t<C^*<t-C_1$ as the other cases are similar. Table 1 outlines the possible values that $\boldsymbol{\psi}_{\mathrm{SU}}$, $\boldsymbol{\psi}^*$ and $(\boldsymbol{\psi}_{\mathrm{SU}}-\boldsymbol{\psi}^*)(\mathbf{1}+(b+1)\mathbf{v})$ can take on for the possible values of $z_1$, $z_2=t-z_1$ and $\mathbf{v}$. Also Figure 2 is helpful. We let $W(\mathbf{z};\mathbf{v})=(\boldsymbol{\psi}_{\mathrm{SU}}(\mathbf{z})-\boldsymbol{\psi}^*(\mathbf{z}))'(\mathbf{1}-(b+1)\mathbf{v})$ and study $E_{\boldsymbol{\mu}}\{W(\mathbf{Z};\mathbf{v})|Z_1+Z_2=t\}$ as a function of $\boldsymbol{\mu}$ for each $\mathbf{v}\in\Gamma$. Note that $\boldsymbol{\mu}$ is in the parameter space only when $\boldsymbol{\mu}\in\Omega_{\mathbf{v}}$.

Using the values from Table 1, it is easy to check that the definition of $C^*$ implies $E_{\mu_1=\mu_2}\{W(\mathbf{Z};\mathbf{v})|Z_1+Z_2=t\}=0$, all $\mathbf{v}\in\Gamma$. For example, say $\mathbf{v}=(0,1)'$. Then, as $Z_1|Z_1+Z_2=t\sim N(\frac{t}{2},\frac{(1-\rho)}{2})$ when $\mu_1=\mu_2$,

$$
\begin{aligned}
&E_{\mu_1=\mu_2}\{W(\mathbf{Z};\mathbf{v})|Z_1+Z_2=t\}\\
&\quad=b\left[\Phi\left(\frac{2C_1-t}{\sqrt{2(1-\rho)}}\right)-\Phi\left(\frac{t-2C_2}{\sqrt{2(1-\rho)}}\right)\right]\\
&\qquad+\left[\Phi\left(\frac{t-2C^*}{\sqrt{2(1-\rho)}}\right)-\Phi\left(\frac{2C_1-t}{\sqrt{2(1-\rho)}}\right)\right]\\
&\qquad+(1-b)\left[\Phi\left(\frac{2C^*-t}{\sqrt{2(1-\rho)}}\right)-\Phi\left(\frac{t-2C^*}{\sqrt{2(1-\rho)}}\right)\right]\\
&\qquad-\left[\Phi\left(\frac{t-2C_1}{\sqrt{2(1-\rho)}}\right)-\Phi\left(\frac{2C^*-t}{\sqrt{2(1-\rho)}}\right)\right]\\
&\qquad-\left[\Phi\left(\frac{2C_2-t}{\sqrt{2(1-\rho)}}\right)-\Phi\left(\frac{t-2C_1}{\sqrt{2(1-\rho)}}\right)\right]
\end{aligned}
$$

$$= 0$$

as $C^*$ is defined by (5.3).

When $\mu_1 \neq \mu_2$ the conditional distribution of $Z_1|Z_1 + Z_2 = t$ is $N(\frac{t}{2} + \eta, \frac{(1-\rho)}{2})$, where $\eta = \mu_1 - \mu_2$.

We further note that $\eta < 0$ when $\boldsymbol{\mu} \in \Gamma_{(0,1)}$ and $\eta > 0$ when $\boldsymbol{\mu} \in \Gamma_{(1,0)}$. The proof can be completed by studying the pattern of sign changes (see Table 1) of $W((z_1, t - z_1); \mathbf{v})$ as a function of $z_1$. It follows from the variation diminishing property [Brown, Johnstone and MacGibbon (1981)] of the normal distribution that $E_{\boldsymbol{\mu}}\{W(\mathbf{Z}; \mathbf{v})|Z_1 + Z_2 = t\} \geq 0$ for all $\boldsymbol{\mu} \in \Omega_{\mathbf{v}}$, $\mathbf{v} \in \Gamma$. This completes the proof. □

**Acknowledgment.** The authors are grateful to the referee who suggested we derive our results for the vector loss function of (2.6). The significance of our results has been enhanced by his suggestion.

Department of Statistics
Rutgers— the State University
of New Jersey
110 Frelinghuysen Road
Piscataway, New Jersey 08854
USA
e-mail: artcohen@rci.rutgers.edu
e-mail: sackrowi@rci.rutgers.edu